\documentclass[conference]{IEEEtran}
\usepackage{ifpdf}

\usepackage{cite}

\ifCLASSINFOpdf
   \usepackage[pdftex]{graphicx}
   
\else
\fi
\usepackage[cmex10]{amsmath}
\usepackage{amssymb}
\usepackage{algorithmic}

\usepackage{array}
\usepackage{cases}

\usepackage{mdwmath}
\usepackage{mdwtab}

\usepackage[tight,footnotesize]{subfigure}

\usepackage{comment}

\begin{document}
%
\title{Examining the rank of Semidefinite Programming for Power System State Estimation}


\author{\IEEEauthorblockN{Byungkwon Park}
\IEEEauthorblockA{Department of Electrical and Computer Engineering\\
University of Wisconsin-Madison\\
Email: bpark52@wisc.edu} \\
}


%


\maketitle

\begin{abstract}
In the power system, state estimation (SE) is important monitoring task for the reliable operation of the system. The optimal estimate from the SE is delivered to all EMS application such as fault analysis, automatic generation control. Hence, it is crucial to have “good” estimation before taking any critical actions. However, the SE problem is challenging problem due to nonconvexity of power flow equations in the nonlinear AC power flow model, which give us a usually local solution. To deal with this nonconvexity, some recent literatures applied the convex semidefinite (SDP) relaxation technique to relax the SE problem attaining or approximating a global solution. In this paper, we investigate the rank of this technique, which is critical to yield a physically meaningful solution with the five-bus test system and propose new approach to possibly reduce the rank by complementing the traditional set of measurement with PMU data. Numerical tests on the standard IEEE 14, 30, 57, 118-bus test system are presented for the demonstration.
\end{abstract}


%
\IEEEpeerreviewmaketitle

\section{Introduction}
The state estimation (SE) play a key role in supervisory control and planning of the power systems. It serves to monitor the state of the power system and estimates the power system states as close as possible to the actual operating state. Therefore, the power system state estimation is the important basis of energy management systems (EMS) for advanced applications in the power system operation and control such as the economic dispatch (ED), automatic generation control (AGC), security analysis (SC), and etc.\cite{MONTICELLI2002}.

Generally, to estimate the power system state variables, complex bus voltages at all buses, the weighted least square (WSL) is employed using the Newton’s Raphson method, which requires the iterative procedure to find a solution \cite{Grainger1994}. Due to its nonconvexity of power flow model relating measurements to state variables, the power system SE problem is inherently nonconvex optimization problem, which provides no guarantee of obtaining a global solution. In the view of dealing with a local minimum, novel approach called relaxation-based semidefinite programming (SDP) has been proposed and applied successfully to various area including the optimal power flow (OPF) problem \cite{LavaeiLow2012} to potentisally obtain a meaningful global solution. Recently, this approach is extended to the SE problem [semidefinite programming for SE, estimating using SDP]. The SDP solution is used to improve Newton’s method by providing a near-global initial guess \cite{WengLiNegiEtAl2012}. Since the WLS formulation is sensitive to bad data, the robust SE problem to handle this sensitivity issue is proposed with the SDP formulation \cite{ZhuGiannakis2012}.

In respect of the ongoing development for the future grid, two major aspects (combination of PMU devices and distributed SE) that will directly impact the SE problem are addressed \cite{HuangWernerHuangEtAl2012}. State estimation model with PMU is proposed and its effect on traditional power system state estimation is analyzed \cite{ChenHanPanEtAl2008}. Compared to traditional measurements, PMU allow us to use linear measurement equations. Modification to account for this linear relationship in the SDP formulation is proposed in \cite{ZhuGiannakis2014}. Furthermore, new regulation and market pricing competition motivate distributed SE to implement interconnection-wide coordinated monitoring. Corresponding algorithms to facilitate distributed SE are presented for the SDP formulation \cite{ZhuGiannakis2012a}, \cite{WengLiNegiEtAl2013}.

Most of approaches above, however, is crucial to have a rank-1 or near rank-1 solution, and it should be questionable and cautious when non-physically meaningful solution \cite{LesieutreMolzahnBordenEtAl2011} could be yielded with a higher rank solution. Therefore, this paper offers insights on the rank of SDP-based SE problem by considering different scenarios. To this end, two factors, the number of measurements and their locations, are considered, and different rank solutions are depicted and discussed with the five-bus test system. Dependency to the number of measurements would be less important than the location of measurement since in practical enough measurements are required to satisfy measurement redundancy not only for detecting and identifying bad data but also for observability \cite{TOYOSHIMACASTILLOFANTINEtAl2012}. However, careful consideration on locating measurements over the power systems would be critical since the significant impact of the location of measurements is shown in this paper. In addition, as deployment of PMUS with the GPS system are accelerating for the future grid, this paper leverage this technique to make the SDP formulation more useful with the attempt to reduce the rank of the SDP-based SE solution. 

The rest of the paper is organized as follows. Section II describes the classical model of the SE problem and Section III constructs the relaxation-based SDP formulation corresponding to the WLS SE problem. Section IV investigate the rank of W matrix with the different number and location of measurements. Also, approach to reduce the rank of W matrix is proposed with the PMU data and conclusion is presented on Section V.


\section{Model/Problem Statement}

Consider the $n$-bus power system with $\mathcal{N}$ denoting the set of all buses, $\mathcal{G}$ denoting the set of all generators, $\mathcal{E}$ representing the set of transmission lines with  $\mathcal{E} \subseteq \mathcal{N}\times\mathcal{N}$. The main purpose of state estimation is to compute the power system state, complex bus voltages $V_ie^{j\delta_i}$ at all buses from a set of redundant measurements $\mathcal{M}$, which are active (reactive) power injection $P_i(Q_i)$, active (reactive) line flow $P_{ij}(Q_{ij})$, bus voltage magnitude $V_i$. 

For AC power system model, the measurement equations are nonlinear and iterative solutions are required as in the Newton-Raphson power-flow procedure. With the power system states $v=[V_ie^{j\delta_i},...,V_ne^{j\delta_n}]^T$, this nonlinear relationship can be captured as \cite{ZimmermanMurillo-ScuteanchezThomas2011}
\begin{subequations}
\label{eq:powerinjection}
\begin{align}
P_{i} &= real[v.*(Y_{bus}v)^*]_i \\
Q_{i} &= imag[v.*(Y_{bus}v)^*]_i \quad i \in \mathcal{N} \\
P_{ij} &= real[v_f.*(Y_{f}v)^*]_{ij} \\
Q_{ij} &= imag[v_f.*(Y_{f}v)^*]_{ij} \quad ij \in \mathcal{E}
\end{align}
\end{subequations}

where $Y_{bus}$ is admittance matrix for bus injections, $Y_f$ is admittance matrix for line power flows, and $v_f$ is corresponding \textquotedblleft from\textquotedblright \, buses for that line. Notice the operator $.*$ indicates the element-wise multiplication in MATLAB. To construct the state estimation problem, those measurements are collected in the vector $z = [z_1,...,z_m]^T$, and the measurement model is formulated as \cite{Grainger1994}
\begin{equation}
\label{eq:model}
z_i = h_i(v)+e_i \quad i \in \mathcal{M}
\end{equation}

where $h_i(v)$ is the nonlinear model described in \eqref{eq:powerinjection} for the measurement of power injections and line flows and linear model for the measurement of voltage magnitudes. It is assumed that $e_i$ is the independent zero-mean normally distributed measurement error with variance $\sigma_i^2$ at the meter $i$. In practical state estimation, the number of actual measurements is far greater than the number of variables required by planning type power flow problem. Therefore, there are many more equations to be solved than unknown state variables i.e., $m \geq n$, and this redundancy is necessary since 1) to guarantee the observability of the whole system 2) there might be bad data or unavailable data due to malfunction in the data-measuring system.

Given $z$, state estimation looks for the best estimated $\hat{x}$ of $x$ according to \eqref{eq:model}. General formulation (weighted least square) for the state estimation is constructed as
\begin{equation}
\label{eq:WLS}
J(x) = \min_{x} \sum_i^m \frac{(z_i-h_i(x))^2}{\sigma_i^2} \quad i \in \mathcal{M} \,\, \text{for} \,\,  x=[\delta_1,...,\delta_n,V_1,...,V_n]
\end{equation}
where $m$ is number of measurements. This problem is denoted as the WLS state estimation and different criteria for the cost function such as absolute value could be used, which is the WLAV estimation to deal with bad data. 

\section{Semidefinite Programming approach for the power system state estimation}

The SE problem is challenging problem due to the nonconvexity of power flow equations. The iterative approach such as Newton’s Raphson is commonly used to solve the problem. This technique has no guarantee of obtaining a global optimum (or even diverge), and is sensitive to an initial guess, and the recent penetration of distributed energy resources makes this sensitivity worse by increasing dynamics of the power systems \cite{GL2014a}.

As illustrated in recent literature on semidefinite programming approaches \cite{LavaeiLow2012} for the standard OPF, coordinate changes applied to problem variables and equations can have significant impact on the efficiency of algorithms and solvers. Not surprisingly, a well-chosen formulation can help in finding better quality solutions, more efficiently. To this end, this section will take advantage of the convex relaxation technique to perform the SE problem. Using this relaxation-based semidefinite programming (SDP) technique, Attempt to attain or approximate a global solution at poly-nominal time complexity is proposed \cite{ZhuGiannakis2011}, \cite{WengLiNegiEtAl2012}. Initial sensitivity described above is handled using the SDP formulation in \cite{WengLiNegiEtAl2012} and Multi-area state estimation with distributed SDP is proposed in [Multi-area]. 

This section formulate the convex relaxation technique to solve the WLS-SE problem \eqref{eq:WLS}. Here we redefine the SDP-related admittance matrices in \cite{LavaeiLow2012} for readers. With basis vector in $R^n$ as $\{e_1,...,e_n\}$ for every $k \in N$ and $(i,j) \in \mathcal{E}$
\begin{align*}
Y_{k} &= e_ke_k^TY_{bus} \\
Y_{ij} &= (\bar{y_{ij}} + y_{ij})e_ie_i^T - y_{ij}e_je_j^T
\end{align*}
related admittance matrices for power injections and line power flows are defined as
\begin{align*}
H_{k,r} &\triangleq \frac{1}{2} \begin{bmatrix} real(Y_k + Y_k^T) & imag(Y_k^T - Y_k) \\ imag(Y_k - Y_k^T) & real(Y_k + Y_k^T) \\ \end{bmatrix} \\
H_{ij,r} &\triangleq \frac{1}{2} \begin{bmatrix} real(Y_{ij} + Y_{ij}^T) & imag(Y_{ij}^T - Y_{ij}) \\ imag(Y_{ij} - Y_{ij}^T) & real(Y_{ij} + Y_{ij}^T) \\ \end{bmatrix} \\
H_{k,q} &\triangleq \frac{-1}{2} \begin{bmatrix} imag(Y_k + Y_k^T) & real(Y_k - Y_k^T) \\ real(Y_k^T - Y_k) & imag(Y_k + Y_k^T) \\ \end{bmatrix} \\
H_{ij,r} &\triangleq \frac{-1}{2} \begin{bmatrix} imag(Y_{ij} + Y_{ij}^T) & real(Y_{ij}^T - Y_{ij}) \\ real(Y_{ij} - Y_{ij}^T) & imag(Y_{ij} + Y_{ij}^T) \\ \end{bmatrix} \\
H_{k,vm} &\triangleq \begin{bmatrix} e_ke_k^T & 0 \\ 0 & e_ke_k^T\\ \end{bmatrix} \\
M_0 &\triangleq \begin{bmatrix} 0 & & & & \cdots & & & 0 \\  & & \ddots & & & & & \\  & & & 0 & & & \\ \colon & & & & 1 & & & \colon \\ & & & & & 0 & & \\ & & & & & & \ddots & & \\ 0 & & & & \cdots & & & 0 \end{bmatrix} \\
X &\triangleq \begin{bmatrix} real(v) \\ imag(v) \\ \end{bmatrix} 
\end{align*}
where $\bar{y_{ij}}$ denotes the value of the shunt element at bus $i$ associated with the line $(i,j)$ and $y_{ij}$ denoting the line admittance between buses $i$ and $j$. Note that $M_0 \in \{0,1\}^{n\times n}$ such that $M_0(i,i) = 1$ if and only if $i$ is correspond to \textquotedblleft n+slack bus\textquotedblright. Using these defined matrices, the linear relationship between measurements z and model with the outer-product matrix $W \triangleq  XX^T$  can be constructed as
\begin{subequations}
\label{eq:SDPMODEL}
\begin{align}
P_{k} &= Tr(H_{k,r}W)   \\
Q_{k} &= Tr(H_{k,q}W)   \\
P_{ij} &= Tr(H_{ij,r}W) \\
Q_{ij} &= Tr(H_{ij,q}W) \\
V_k^2 &= Tr(H_{k,vm}W) 
\end{align}
\end{subequations}
where $tr$ indeicates the trace operator in MATLAB. Then, corresponding optimization problem can be written as
\begin{align}
\label{eq:SDP1}
J(x) = \min_{W} \sum_i^m \frac{(z_i-tr(H_iW))^2}{\sigma_i^2} \quad i \in \mathcal{M} \\
\text{suject to} \,\,  W\succeq 0, rank(W) = 1, tr(M_0W)=0 \nonumber
\end{align}

where $H_i$ is the corresponding matrix for each measurement specified in \eqref{eq:SDPMODEL}. The constraint $tr(M_0W)=0$ is for angle reference condition (zero angle for the slack bus). The positive semi-definiteness and rank constraints is imposed to ensure that there exists a vector  such that $W=XX^T$. The above formulation is not convex due to the objective function being of degree 4 with respect to the variable X and rank constraint. However, by defining some auxiliary variable one can reformulate the problem in a quadratic way with respect to the variable X. To this end, introduce a scalar variable $\alpha_i$ for $i\in \mathcal{M}$ with $\frac{(z_i-tr(Y_iW))^2}{\sigma_i^2} = f_i(tr(H_iW)) \leq \alpha_i$. Then this is equivalent to 
\begin{align*}
\frac{z_i^2 - 2z_itr(H_iW) + tr(H_iW)^2}{\sigma_i^2} - \alpha_i & \leq 0 \\
-(\frac{z_i^2 - 2z_itr(H_iW) + tr(H_iW)^2}{\sigma_i^2}) + \alpha_i & \geq 0 
\end{align*}
This equation can be reformulated by Schur\textquoteright s complement formula as 
\begin{align*}
\begin{bmatrix} \alpha_i - \frac{z_i^2}{\sigma_i^2} + \frac{2z_itr(H_iW)}{\sigma_i^2} & -
\frac{1}{\sigma}tr(H_iW) \\ -
\frac{1}{\sigma}tr(H_iW) & 1 \\ \end{bmatrix} \succeq 0 \quad i \in \mathcal{M}
\end{align*}

Thus, final equivalent optimization problem is 
\begin{equation}
\label{eq:SDP1}
\min_{W,\alpha_i} \sum_i^m \alpha_i 
\end{equation}
subject to
\begin{align}
\begin{bmatrix} \alpha_i - \frac{z_i^2}{\sigma_i^2} + \frac{2z_itr(H_iW)}{\sigma_i^2} & -
\frac{1}{\sigma}tr(H_iW) \\ -
\frac{1}{\sigma}tr(H_iW) & 1 \\ \end{bmatrix} \succeq 0 \quad i \in \mathcal{M} \nonumber  \\
W\succeq 0, rank(W) = 1, tr(M_0W)=0 \nonumber
\end{align}

Notice that since optimization problem \eqref{eq:SDP1} has a rank constraint, it is still nonconvex. However, by removing the rank constraint, we obtain a semidefinite programming, which leads us to have the following convex optimization problem.
\begin{equation}
\label{eq:SDP2}
\min_{W,\alpha_i} \sum_i^m \alpha_i 
\end{equation}
subject to
\begin{align}
\begin{bmatrix} \alpha_i - \frac{z_i^2}{\sigma_i^2} + \frac{2z_itr(H_iW)}{\sigma_i^2} & -
\frac{1}{\sigma}tr(H_iW) \\ -
\frac{1}{\sigma}tr(H_iW) & 1 \\ \end{bmatrix} \succeq 0 \quad i \in \mathcal{M} \nonumber  \\
W\succeq 0, tr(M_0W)=0 \nonumber
\end{align}

This optimization problem \eqref{eq:SDP2} is a relaxation of the optimization problem \eqref{eq:SDP1} and does not always have a rank-1 solution, which is no longer physically meaningful \cite{LesieutreMolzahnBordenEtAl2011}. We investigate this rank condition of W matrix for the SE problem with the five-bus test system in the next chapter.  

\section{Numerical Example}

Here we provide a simple computational example to analyze the case whether rank-1 solution is obtained or not. The SDP formulation of the SE problem was solved using MATLAB-based optimization package CVX \cite{GrantBoyd2011}. 

\subsection{Methodology}
In this section, numerical tests are implemented on the five-bus system depicted in Fig. 1. The set of measurements are collected first by solving the power flow solutions that provide the noise-free measurements (true measurements), and then randomly extracts (large enough to ensure the system observability) those values to have a set of measurement $z$. These measurement values from power flow solution are shown in Tables I and II. All measurements are corrupted by Independent Gaussian noise with  $\sigma_i$ equal to 0.008 at line power flow meters, 0.01 at injection power meters, and 0.004 at voltage meters.

\begin{table} 
	[htb]
	\renewcommand{\arraystretch}{1.1}
	\centering
	\begin{tabular}{|c|c|c|}
		\hline
		& Power Injection & Bus Voltage Phasor \\
		\hline\hline
		Bus 1 & 4.5835 + j0.9233 & 1.0000 + j0.0000 \\
		\hline
		Bus 2 & -1.3000 - j0.2000 & 0.8589 - j0.1930 \\
		\hline
		Bus 3 & -1.3000 - j0.2000 & 0.8485 - j0.2006 \\
		\hline
		Bus 4 & -0.6500 - j0.1000 & 0.2863 - j0.8692 \\
		\hline
		bus 5 & 0.0000 + j1.2371 & 0.3069 - j0.9518 \\
		\hline
		\noalign{\vskip 2mm}   
	\end{tabular}
	\caption{Measurement Values from Power Flow Solution}
\end{table}
\vspace{-0.5cm}
\begin{table} 
	[htb]
	\renewcommand{\arraystretch}{1.1}
	\centering
	\begin{tabular}{|c|c|c|}
		\hline
		From bus & To bus & Line Power Flow \\
		\hline\hline
		Bus 1 & Bus 2  & 2.3729 + j0.5131 \\
		\hline
		Bus 1 & Bus 3  & 2.2106 + j0.4102 \\
		\hline
		Bus 2 & Bus 4  & 0.7379 - j0.2239 \\
		\hline
		Bus 3 & Bus 5  & 0.7562 - j0.2900 \\
		\hline
		Bus 4 & Bus 5  & -0.3003 - j0.5962 \\
		\hline
		Bus 2 & Bus 3  & 0.0992 + j0.0065 \\
		\hline
		\noalign{\vskip 2mm}   
	\end{tabular}
	\caption{Measurement Values from Power Flow Solution}
\end{table}

\begin{figure}[htb]
	\centering
	\includegraphics[height=1.5in,width=2.5in]{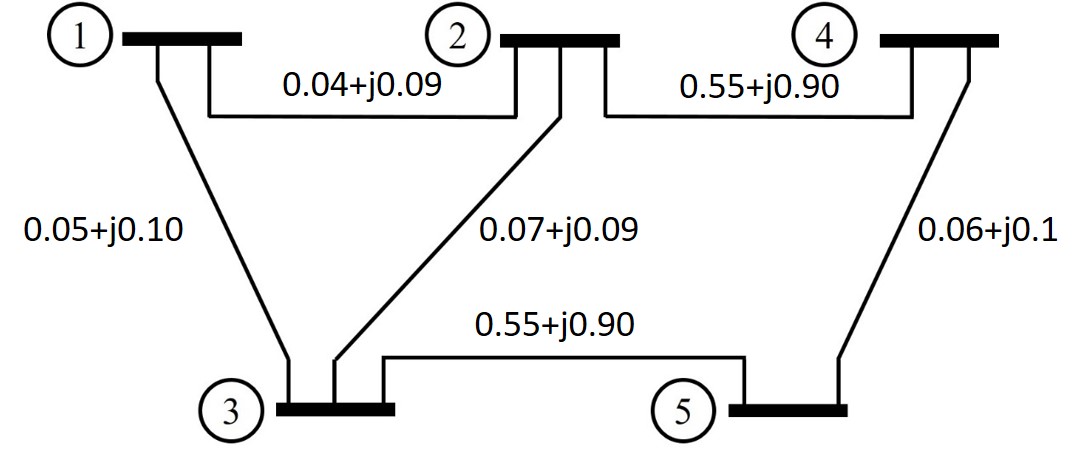}
	\caption{Five-bus test system}
	\label{fig:bus-compare}
\end{figure}

\subsection{Number of measurement vs the Rank of W matrix}
We first discuss number of measurements affecting the rank of W matrix in the five-bus system. Number of measurements play an important role to obtain a physically meaningful global solution, and one of the assumptions discussed in \cite{ZhuGiannakis2014} is that every bus is equipped with a voltage magnitude meter to satisfy rank(W)=1. Although bus voltage magnitude meters are deployed in most substations to monitor voltage rating for the system stability, this condition is approximation of the realistic power system. Instead, we increase the number of measurements from 40\% to 70\% for each quantity and investigate the rank of W matrix. 

\begin{figure}[htb]
	\centering
	\includegraphics[height=1.5in,width=3in]{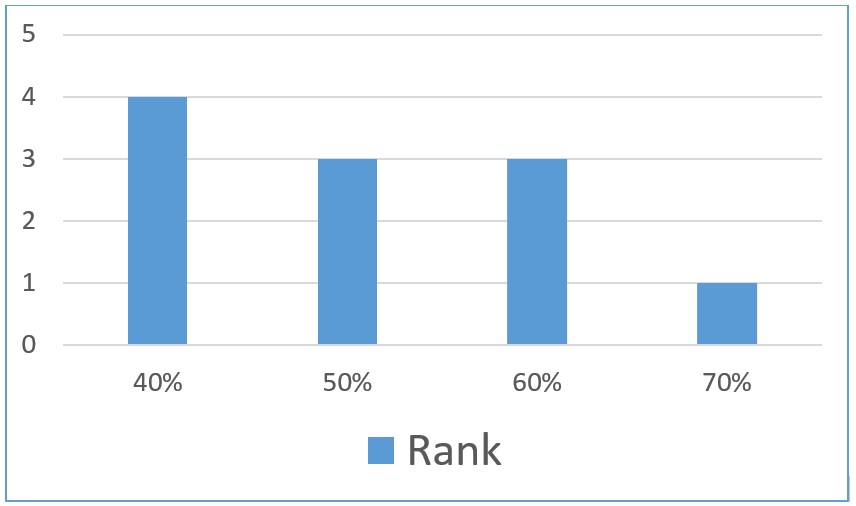}
	\caption{The rank of W matrix with measurements increased from 40\% to 70\%}
	\label{fig:bus-compare}
\end{figure}

In Fig 2, it clearly describes that more measurements tend to obtain a solution with a lower rank W (hence near-global solution to the WLS SE problem). The SDP formulation provides non-meaningful solutions until 60\% with higher rank of W matrix than 1, and yield a rank-1 solution with 60\% and 70\%. It appears that measurement redundancy is important for the SDP formulation in terms of the rank of W matrix as well as observability and bad data.

\subsection{Locational dependency of the Rank of W matrix}
We now discuss the second factors, location of measurements, affecting the rank of W matrix in the five-bus system. Change the location of measurements while leaving the number of measurements unchanged (consist of 70\% values for each quantity form Tables I and II). Specific locations of measurements are shown in Fig. 3(a) and Fig. 3(b). Location corresponding to rank-1 solution (physically meaningful) is shown in Fig. 3(a), whereas non-physically meaningful solution is shown in Fig. 3(b). Only difference between two locations is the bus power injection at bus 3 so that only one quantity (line power flow) is measured around at bus 3 in rank-3 solution. Notice that bus 3 is the critical location as it is connected with 3 lines implying many quantities are associated with, and unrealistically large value is discovered at bus 3 for rank-3 solution. Solution corresponding to these locations with the SDP and WLS SE problem is described in Tables III and IV respectively.

\begin{figure}[!htb]
	\centering
	\subfigure[Rank-1 solution]{%
		\includegraphics[height=1in,width=0.45\linewidth]{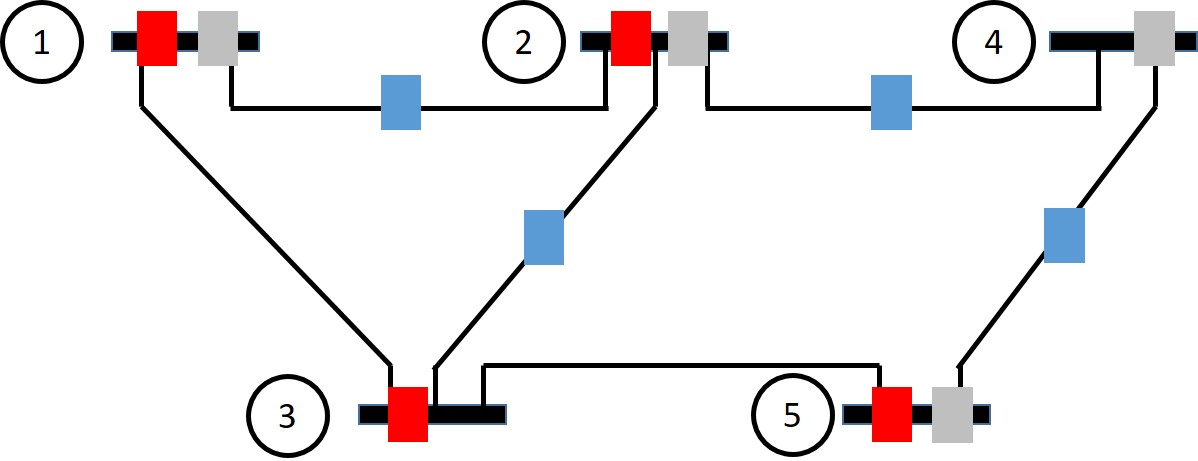}
		\label{fig:subfigure1}}
	\quad
	\subfigure[Rank-3 solution]{%
		\includegraphics[height=1in,width=0.45\linewidth]{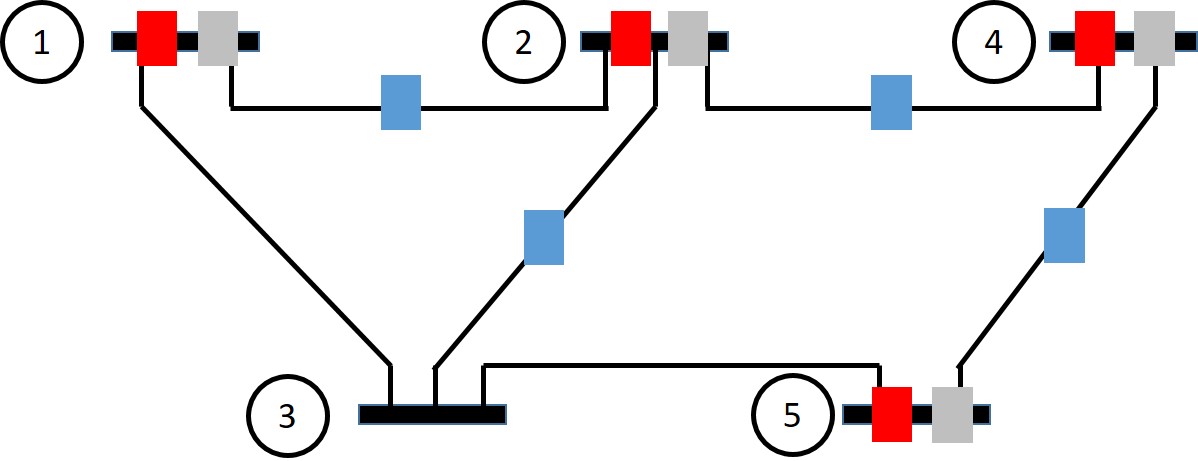}
		\label{fig:subfigure2}}
	\caption{Five-bus system. Red, gray and blue squares indicate the power injections, voltage magnitude and line power flow meter locations respectively}
	\label{fig:figure}
\end{figure}

\begin{table} 
	[htb]
	\renewcommand{\arraystretch}{1.1}
	\centering
	\begin{tabular}{|c|c|c|}
		\hline
		 & SDP (rank-1) & SDP (rank-3) \\
		\hline\hline
		Bus 1 & 0.9990 + j0.0000 & 0.0047 - j0.0000 \\
		\hline
		Bus 2 & 0.8578 - j0.1930 & 0.0041 - j0.0009 \\
		\hline
		Bus 3 & 0.8477 - j0.2005 & 46.804 - j102.59 \\
		\hline
		Bus 4 & 0.2932 - j0.8687 & 0.0014 - j0.0041 \\
		\hline
		bus 5 & 0.3138 - j0.9504 & 0.0015 - j0.0045 \\
		\hline
		Sum square error & 8.3105 & 6.5951 \\
		\hline
		\noalign{\vskip 2mm}   
	\end{tabular}
	\caption{Solution to Five-bus test system wih the SDP fomulation}
\end{table}

\begin{table} 
	[!htb]
	\renewcommand{\arraystretch}{1.1}
	\centering
	\begin{tabular}{|c|c|c|}
		\hline
		& WLS (rank-1) & WLS (rank-3) \\
		\hline\hline
		Bus 1 & 0.9985 + j0.0000 & 0.9989 + j0.0000 \\
		\hline
		Bus 2 & 0.8573 - j0.1932 & 0.8578 - j0.1931 \\
		\hline
		Bus 3 & 0.8471 - j0.2006 & 0.8477 - j0.2006 \\
		\hline
		Bus 4 & 0.2933 - j0.8699 & 0.2939 - j0.8695 \\
		\hline
		bus 5 & 0.3138 - j0.9515 & 0.3147 - j0.9511 \\
		\hline
		Sum square error & 8.4447 & 8.8790 \\
		\hline
		\noalign{\vskip 2mm} 
	\end{tabular}
	\caption{Solution to Five-bus test system wih the WLS fomulation}
\end{table}

Table III shows solutions from the SDP formulation and Table IV shows solutions from the WLS formulation. Notice that the solution (as best approximated rank-1 solution) in Table III is recovered by eigenvalue decomposition using the largest eigenvalue such as $v=\sqrt{\lambda_1}u_1$ with $\lambda_1$ is the largest eigenvalue and $u_1$ is the corresponding eigenvector of W matrix. In Table III, first column yields a physically meaningful result as evidenced by rank-1 and non-physically meaningful solution in second column with rank-3, whereas the WLS formulations still provide a physically meaningful solution regardless of the location of measurements in this example.

The optimal objective value for each formulation is shown the last row, which demonstrates the SDP solution is less than the WLS solution. Thus, the objective function values at the SDP formulation lower bounds that of the WLS formulation (classical formulation) as expected. Other test system including 14-bus, 30-bus, 57-bus and 118-bus systems are examined with fixed number of measurement in Fig 4. In this work, 70\% of power flow solution are collected to ensure the observability of the power systems. For the fair analysis, problems are solved 50 times and average value is taken to discuss the rank of W matrix.

\begin{figure}[!htb]
	\centering
	\includegraphics[height=2in,width=3in]{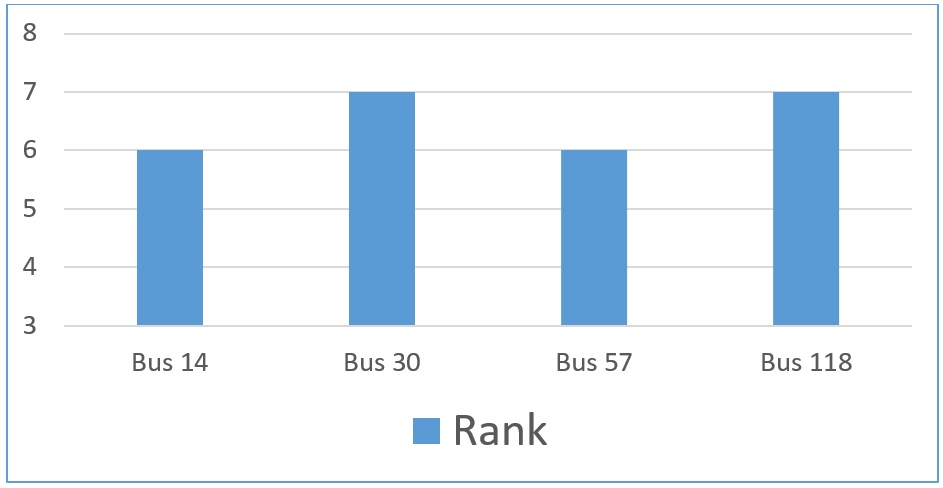}
	\caption{The rank of W matrix with different test systems (average value)}
	\label{fig:bus-compare}
\end{figure}

As depicted in Fig 4, none of the systems produces the rank-1 solution, but still useful information can be considered, particularly in cases for which W matrix is close to a rank-1 matrix. Closest rank one approximation to W matrix could be served as an initial condition to the WLS SE problem.  



\subsection{Reducing the Rank of W matrix}

Here we attempt to reduce the rank of W matrix to be as close as to a rank-1 matrix. As mentioned in the previous section, the rank of W matrix is also dependent on the number of measurement. With the introduction of PMU, now the power system states can be measured directly with increased precision in measuring phasor angles due to GPS’s network-wide synchronization and short measurement periodic time. This suggests that we could complement the set of measurements with PMU data to perform the SE problem. To this end, angle equation with the rectangular formulation for each bus is modified as
\begin{align}
\delta_i & = tan^{-1}\Big(\frac{V_i^q}{V_i^d}\Big) \nonumber \\
tan(\delta_i) & = tan\Big(tan^{-1}\Big(\frac{V_i^q}{V_i^d}\Big)\Big) = \frac{V_i^q}{V_i^d} \nonumber \\
(V_i^d)^2tan(\delta_i) & = V_i^qV_i^d \nonumber \\
(V_i^d)^2tan(\delta_i) & - V_i^qV_i^d  = 0 \nonumber \\
\text{this is equivalent to} \nonumber \\
tr(H_{i,a}W) & = 0 \,\,\, i \in \{\text{buses with angle measurement}\}
\label{eq:ANGLE}
\end{align}

Note that $H_{i,a} \in \{tan(\delta_i),-1,0\}^{2n\times 2n}$ such that $H_{i,a}(i,i) = tan(\delta_i)$, $H_{i,a}(i,i+n) = -1$, and 0 otherwise. Our purpose here is to have the equation relating angle measurements to state variables in the quadratic way. To this end, each angle measurement is expressed as function of state variables in the fraction with arctangent formula. By taking the tangent, the inverse of arctangent function, variables in the fraction are pulled out. Last step is to multiply both sides by $V_d^2$ to be expressed in the quadratic way.

Angle equation including variables in the fraction is converted to the one including variables in the quadratic way, which fits the SDP formulation. It is assumed that voltage angle is noise-free measurement with high precision from PMU. The equation (8) is imposed to only the buses equipped with bus voltage magnitude measurement in the original set of measurements. To see the impact on the rank, we gradually increase number of angle measurements from 0\% to 30\%. For the fair analysis, problems are solved 50 times and the lowest rank value taken to discuss since we are interested in the best scenario that can be improved.

\begin{figure}[htb]
	\centering
	\includegraphics[height=2in,width=3in]{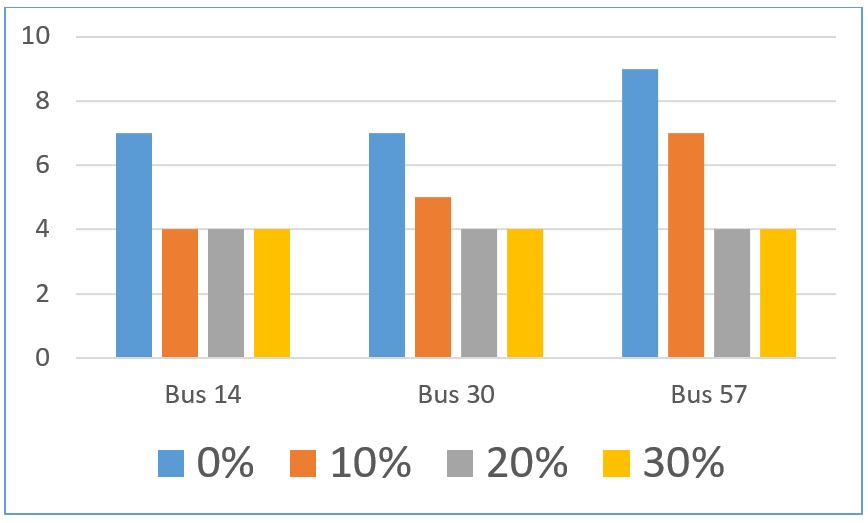}
	\caption{Two port representation and its matrix form}
	\label{fig:bus-compare}
\end{figure}

The rank of W matrix with angle measurements for each test system is shown in Fig 5. With this approach, some, but not all, of test systems yield a lower rank of W matrix than without angle measurements. This example also demonstrates that more measurements tend to have a lower rank of W matrix. As noted previously, near rank-1 approximated solution could be used as an initial condition for the WLS SE problem to possibly obtain a global solution. We anticipate that this little improvement on the rank shown here could help to have a near rank-1 approximated solution for further applications.

\subsection{Discussion}
Results above establish that the rank of W matrix is dependent on not only the number of measurement but also their locations. In this work, true measurements from power flow solution are corrupted with zero mean Gaussian noise, those noises are fixed for each simulation. In fact, the rank of W matrix even varies with different values of noise. This is a very unattractive behavior since for example one bad data corrupted by noise would change the rank of W matrix, so it could yield a less accruate approximated rank-1 solution. Therefore, bad data detection would be crucial for the implementation of the SDP formulation to the WLS SE problem. 

For the test, randomly generated noises are added to the true measurements (power flow solutions). All (100\%) noise-added measurements are used to check the impact of noise on the rank of W matrix for the 14-bus test system. Test is simulated 50 times and in Fig 6, it is shown that the rank of W matrix changes from 3 as the best case to 6 as the worst case, which implies that some badly noise-added data (we call this as bad data) is needed to be detected and removed for more accurate result of the SDP formulation.  

\begin{figure}[htb]
	\centering
	\includegraphics[height=2in,width=3in]{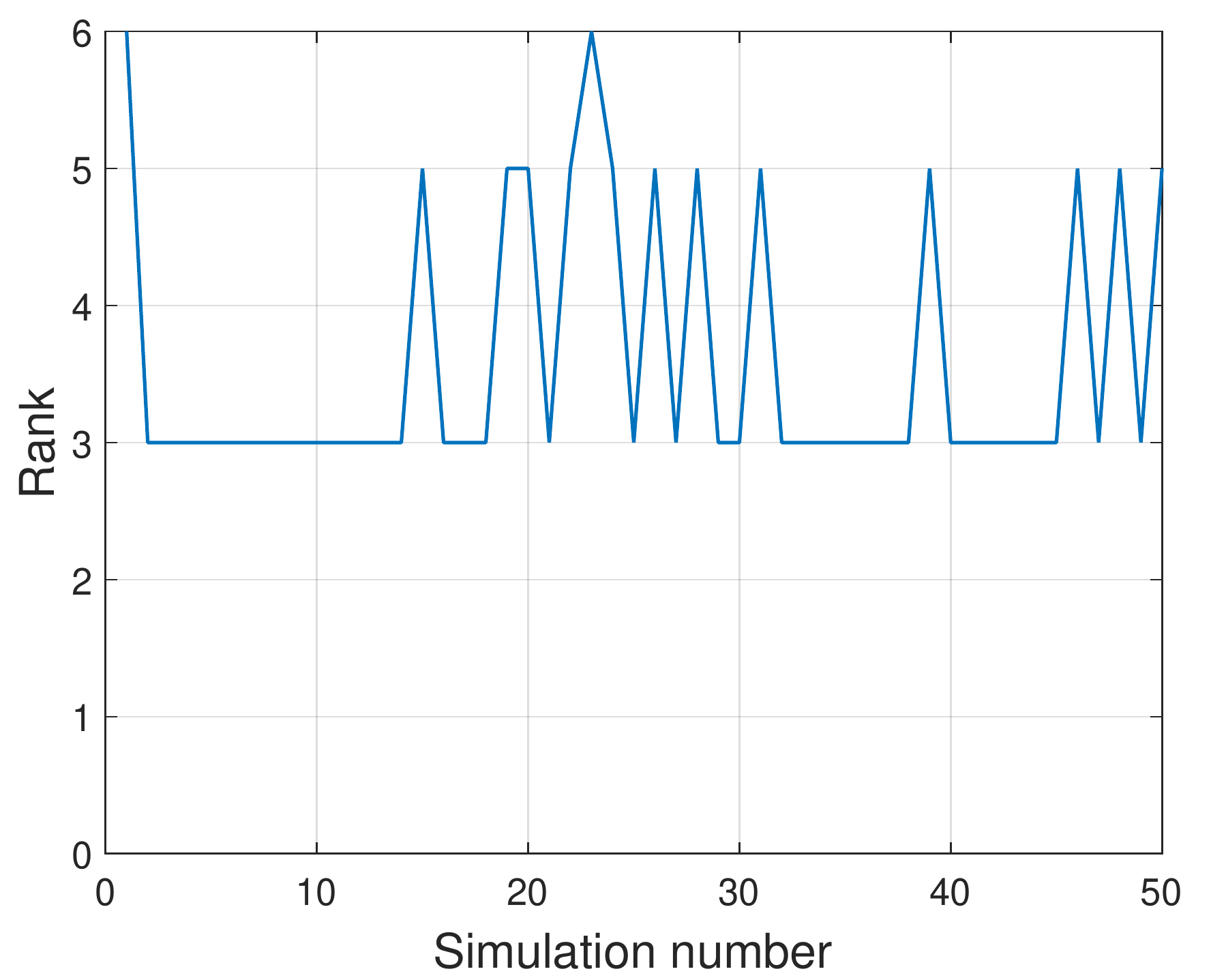}
	\caption{Change of the rank of W matrix with randomly generator noises for the 14-bus test system}
	\label{fig:bus-compare}
\end{figure}

Another issue to be careful is to consider the ratio of the largest eigenvalue to other eigenvalues. Let $W_2$ be a rank-2 matrix with two non-zero eigenvalues $(15,13)$ and $W_5$ be a rank-5 matrix with five non-zero eigenvalues $(29,2,0.5,0.1,0.1)$. In this case, the ratio of the largest eigenvalue to other eigenvalues is more important than the rank of that matrix. Observe the ratio $\frac{\lambda_2}{\lambda_1}=0.8667$ for $W_2$ matrix and $\frac{\sum_j^5\lambda_j}{\lambda_1}=0.091$ for $W_5$ matrix. The best approximated rank-1 solution is recovered with the largest eigenvalue and this would be more accurate when other eigenvalues are relatively very small compared to the largest eigenvalue. Hence, $W_5$ matrix yield a better approximated rank-1 solution, even though its rank is higher than $W_2$ matrix. This situation happens sometime but not often so careful consideration would be necessary when the SDP formulation is used as an initial condition for the WLS SE problem.

\section{Conclusion}
This paper has investigated existing applications of the semidefinite programming for the state estimation problem. First we have constructed the SDP formulation for the WLS SE problem and discussed two factors that might affects the rank of W matrix: the number of measurements and location of measurements. In terms of the location of measurements, this paper has provided an example that the SDP formulation may fail to obtain physically meaningful results with the same number of measurements in the power system.

To explore the possible improvement on the rank of W matrix to have a near rank-1 approximated solution, recent deployment of the PMU data is considered, allowing us to access angle measurements at each bus. Rectangular formulation including variables in the fraction is modified to one that has variables in the quadratic way. Number of angle measurements are then gradually increased from 0\% to 40\% to discuss its impact on the rank of W matrix.  Although not all of test cases was successful in reducing the rank, some test cases able to reduce the rank were described, which gives a better approximated rank-1 solution. Finally, concerns about the rank with bad data and recovering of solution are discussed. 

The SDP formulation in this paper show that it is not capable of finding a physically meaningful solution for some cases. However, the use of relaxation-based SDP methods for the SE problem has significant potentials and further enhancements to the rank of W matrix to have a near rank-1 solution would be the next approach for further applications. Also, the distributed SDP formulation can be applied to reduce the burden of computational issues, which only obtains submatrices of W matrix. This implies that second level optimization to reduce the rank of W matrix could be considered to have a better approximated rank-1 solution after the distributed SDP formulations is solved.




%

\bibliographystyle{IEEEtran}
\bibliography{SDPSE}

\end{document}